\documentstyle[12pt]{report}
\newcommand{\nl}{\mbox{}\\}
\newcommand{\np}{\par \noindent \mbox{\mbox{}\hspace{1.0cm}}}
\newcommand{\intR}{\mbox{$ {\displaystyle \int_{- \infty}^{+ \infty} }$}}
\newcommand{\intline}{\intR}
\newcommand{\tend}{\mbox{$ \rightarrow $}}
\newcommand{\limit}{\mbox{$ {\displaystyle
            \lim_{\mbox{}_{\scriptstyle t \,\rightarrow\, \infty }} }$}}
\newcommand{\limitsup}{\mbox{$ {\displaystyle
            \limsup_{\mbox{}_{\scriptstyle t \,\rightarrow\, \infty }} }$}}
\newcommand{\limitinf}{\mbox{$ {\displaystyle
            \liminf_{\mbox{}_{\scriptstyle  t \,\rightarrow\, \infty }} }$}}
\newcommand{\Lone}{\mbox{$ L^{{\scriptstyle 1}}({\bf R}) $}}

\newcommand{\Linfty}{\mbox{$ L^{{\scriptstyle \infty}}({\bf R}) $}}
\newcommand{\nLone}{\mbox{$ \|_{\mbox{}_{\mbox{}_{\scriptstyle 1}}} $}}
\newcommand{\nLtwo}{\mbox{$ \|_{\mbox{}_{\mbox{}_{\scriptstyle 2}}} $}}
\newcommand{\nLp}{\mbox{$ \|_{\mbox{}_{\mbox{}_{\scriptstyle p}}} $}}
\newcommand{\nLinfty}{\mbox{$ \|_{\mbox{}_{\mbox{}_{\scriptstyle \infty}}}$}}
\newcommand{\nLtwohalf}{\mbox{$ \|_{\mbox{}_{\mbox{}_{\scriptstyle 2}}}
   ^{\mbox{}^{\scriptstyle \frac{\scriptstyle 1}{\scriptstyle 2} }} $}}
\newcommand{\nLoneINTP}{\mbox{$ \|_{\mbox{}_{\mbox{}_{\scriptstyle 1}}}
   ^{\mbox{}^{\scriptstyle \frac{\scriptstyle 1}{\scriptstyle p} }} $}}

\newcommand{\nLpp}{\mbox{$ \|_{\mbox{}_{\mbox{}_{\scriptstyle p}}}
                             ^{\mbox{}^{\scriptstyle p}} $}}
\newcommand{\ut}{\mbox{$ u_{\mbox{}_{\scriptstyle t}} $}}
\newcommand{\ux}{\mbox{$ u_{\mbox{}_{\scriptstyle x}} $}}
\newcommand{\uxx}{\mbox{$ u_{\mbox{}_{\scriptstyle x x}} $}}
\newcommand{\phit}{\mbox{$ {\varphi}_{\mbox{}_{\scriptstyle t}} $}}
\newcommand{\phix}{\mbox{$ {\varphi}_{\mbox{}_{\scriptstyle x}} $}}
\newcommand{\phixx}{\mbox{$ {\varphi}_{\mbox{}_{\scriptstyle x x}} $}}
\newcommand{\phizero}{\mbox{$ {\varphi}_{\mbox{}_{\scriptstyle 0}} $}}

\newcommand{\psix}{\mbox{$ {\psi}_{\mbox{}_{\scriptstyle x}} $}}

\newcommand{\omegat}{\mbox{$ {\omega}_{\mbox{}_{\scriptstyle t}} $}}

\newcommand{\omegaxx}{\mbox{$ {\omega}_{\mbox{}_{\scriptstyle x x}} $}}

\newcommand{\factort}{\mbox{$ {\displaystyle
    \frac{1}{\sqrt{\,4\,\pi\,c\,t\;}\,} }$}}
\newcommand{\factormt}{\mbox{$ {\displaystyle
    \frac{m}{\sqrt{\,4\,\pi\,c\,t\;}\,} }$}}
\newcommand{\factormm}{\mbox{$ {\displaystyle
    \frac{|\,m\,|}{\sqrt{\,4\,\pi\,c\,t\;}\,}  }$}}
\newcommand{\factorpi}{\mbox{$ {\displaystyle
    \frac{1}{\sqrt{\,\pi\;}\,} }$ }}
\newcommand{\factortm}{\mbox{$ {\displaystyle
    \frac{m}{\sqrt{\,4\,\pi\,c\,t\;}\,} }$}}
\newcommand{\expxy}{\mbox{$ {\displaystyle
    \mbox{\sf e}^{\mbox{}^{\scriptstyle - \,
       \mbox{\footnotesize $ {\displaystyle
       \frac{\mbox{$ \scriptstyle \,(\,\mbox{\footnotesize $x$} \,-\,
       \mbox{\footnotesize $y$} \,)^{\mbox{}^{\scriptstyle 2}} $}\!\!\!}
            {\mbox{$ \scriptstyle
       \mbox{\footnotesize $4$} \, \mbox{\footnotesize $c$} \,
       \mbox{\footnotesize $t$} $}} }$} }} }$}}
\newcommand{\expxyat}{\mbox{$ {\displaystyle
    \mbox{\sf e}^{\mbox{}^{\scriptstyle - \,
       \mbox{\footnotesize $ {\displaystyle
       \frac{\mbox{$ \scriptstyle \,(\,\mbox{\footnotesize $x$} \,-\,
       \mbox{\footnotesize $y$} \,-\, \mbox{\footnotesize $a$} \,
       \mbox{\footnotesize $t$}\,)^{\mbox{}^{\scriptstyle 2}} $}\!\!\!}
            {\mbox{$ \scriptstyle \mbox{\footnotesize $4$} \,
       \mbox{\footnotesize $c$} \, \mbox{\footnotesize $t$} $}} }$} }} }$}}
\newcommand{\expxiy}{\mbox{$ {\displaystyle
    \mbox{\sf e}^{\mbox{}^{\scriptstyle - \,
       \left(\, \mbox{\footnotesize $\xi$} \,-\,
            \mbox{\footnotesize $ {\displaystyle
             \frac{\mbox{\footnotesize $y$}}
             {\mbox{$ \scriptstyle \,
              \sqrt{\,\mbox{\footnotesize $4$} \, \mbox{\footnotesize $c$}
                    \,\mbox{\footnotesize $t$} \;}\,$}} }$}
       \right)^{\scriptstyle 2} }} }$ }}
\newcommand{\expxiya}{\mbox{$ {\displaystyle
    \mbox{\sf e}^{\mbox{}^{\scriptstyle - \,
       \left(\, \mbox{\footnotesize $\xi$} \,+\,
          \mbox{\footnotesize $ {\displaystyle
          \frac{\mbox{$ \scriptstyle \, \mbox{\footnotesize $\alpha$} \,-\,
                \mbox{\footnotesize $y$} \,$}}
            {\mbox{$ \scriptstyle \,\sqrt{\, \mbox{\footnotesize $4$} \,
                \mbox{\footnotesize $c$} \, \mbox{\footnotesize $t$} \,}\,$}}
       }$} \right)^{\scriptstyle 2} }} }$ }}
\newcommand{\expxiyan}{\mbox{$ {\displaystyle
    \mbox{\sf e}^{\mbox{}^{\scriptstyle - \,
       \left(\, \mbox{\footnotesize $\xi$}_{\mbox{}_{\scriptstyle n}} \,+\,
         \mbox{\footnotesize $ {\displaystyle
         \frac{\mbox{$ \scriptstyle \,\mbox{\footnotesize $\alpha$} \,-\,
                       \mbox{\footnotesize $y$} \,$}}
              {\mbox{$ \scriptstyle \,\sqrt{\, \mbox{\footnotesize $4$} \,
                       \mbox{\footnotesize $c$} \,
                \mbox{\footnotesize $t$}_{\mbox{}_{\scriptstyle n}}\,}\,$}}
       }$} \right)^{\scriptstyle 2} }} }$ }}
\newcommand{\expxi}{\mbox{$ {\displaystyle
    \mbox{\sf e}^{\mbox{}^{\mbox{$\scriptstyle - \,
       \mbox{\footnotesize $\xi$}^{\mbox{}^{\scriptstyle 2}} $}}} }$ }}
\newcommand{\ebm}{\mbox{$ {\displaystyle
    \mbox{\sf e}^{\mbox{}^{\scriptstyle - \,
      \mbox{\footnotesize $ {\displaystyle
      \frac{ \,\mbox{$ \scriptstyle \mbox{\footnotesize $b$} \,
                                    \mbox{\footnotesize $m$} $} \,}
           { \,\mbox{$ \scriptstyle \mbox{\footnotesize $2$} \,
                                    \mbox{\footnotesize $c$} $} \,}
           }$} }} }$ }}
\newcommand{\fbcm}{\mbox{$ {\displaystyle
    \mbox{\huge $|$}\, \frac{\,2\,c\,}{\,b\,m\,} \:
      \mbox{\Large $($} \, 1 \,-\, \ebm \,\mbox{\Large $)$}
    \; \mbox{\huge $|$} }$ }}
\newcommand{\fourcp}{\mbox{$ {\displaystyle
    (\, 4\,c\,)^{\mbox{}^{\mbox{}^{\scriptstyle \!
    \frac{ \mbox{$ \scriptstyle 1 $} }
         { \mbox{$ \scriptstyle \,2\,p\, $} } }}} }$}}
\newcommand{\tmp}{\mbox{$ {\displaystyle t^{\mbox{}^{\scriptstyle
    \,- \frac{ \mbox{$ \scriptstyle 1 $} }{ \mbox{$ \scriptstyle 2 $} }
      \left( \frac{\mbox{}}{\mbox{}} \! 1 \right.
        \,-\, \frac{ \mbox{$ \scriptstyle 1 $} }
                   { \mbox{$ \scriptstyle p $} }
      \left. \frac{\mbox{}}{\mbox{}} \!\!\right) }} }$}}
\newcommand{\tpp}{\mbox{$ {\displaystyle t^{\mbox{}^{\scriptstyle
      \, \frac{ \mbox{$ \scriptstyle 1 $} }{ \mbox{$ \scriptstyle 2 $} }
      \left( \frac{\mbox{}}{\mbox{}} \! 1 \right.
       \,-\, \frac{ \mbox{$ \scriptstyle 1 $} }
                  { \mbox{$ \scriptstyle p $} }
      \left. \frac{\mbox{}}{\mbox{}} \!\!\right) }} }$}}
\newcommand{\tppp}{\mbox{$ {\displaystyle t^{\mbox{}^{\scriptstyle
      \frac{\scriptstyle p }{\scriptstyle 2 } \,-\,
      \frac{\scriptstyle 1 }{\scriptstyle 2 } }} }$}}
\newcommand{\thalf}{\mbox{$ {\displaystyle t^{\mbox{}^{\scriptstyle
      \frac{\scriptstyle 1 }{\scriptstyle 2 } }} }$}}

\newcommand{\HCxxit}{\mbox{$ {\displaystyle
      {\sf e}^{\mbox{}^{\scriptstyle
         -\, \mbox{\footnotesize $ {\displaystyle
             \frac{\mbox{\footnotesize $b$}}
                 {\scriptstyle \,\mbox{\footnotesize $2$}\,
                                 \mbox{\footnotesize $c$}\,} \,
         \int_{\scriptstyle - \infty}^{\scriptstyle x} \!\!\!
            u(\xi,t) \, \mbox{\footnotesize $d\xi$} }$} }} }$ }}
\newcommand{\HCvxxit}{\mbox{$ {\displaystyle
      {\sf e}^{\mbox{}^{\scriptstyle - \,
         \mbox{\footnotesize $ {\displaystyle
            \frac{\mbox{\footnotesize $b$}}
                 {\,\mbox{\footnotesize $2$}\,\mbox{\footnotesize $c$}\,}
                 \int_{\scriptstyle - \infty}^{\scriptstyle x} \!\!\!
            v(\xi,t) \, \mbox{\footnotesize $d\xi$} }$ } }} }$}}
\newcommand{\HCxxizero}{\mbox{$ {\displaystyle
      {\sf e}^{\mbox{}^{\scriptstyle - \,
          \mbox{\footnotesize $ {\displaystyle
            \frac{\mbox{\footnotesize $b$}}
                 {\,\mbox{\footnotesize $2$}\,\mbox{\footnotesize $c$}\,}
                 \int_{\scriptstyle - \infty}^{\scriptstyle x} \!\!\!
             \mbox{\footnotesize $ u_0(\xi) \, d\xi $} }$ } }} }$}}
\newcommand{\Hzero}{\mbox{$ {\cal H}_{0} $}}
\newcommand{\Hzerox}{\mbox{$ {\cal H}_{0}(x) $}}
\newcommand{\Hzeroy}{\mbox{$ {\cal H}_{0}(y) $}}
\newcommand{\pstar}{\mbox{$ p_{\mbox{}_{\scriptstyle \ast}} $}}

\begin{document}
%
%---------------------------------------------------------------------
%
%                         TITLE / AUTHOR :
%
% \title{\centering \large \bf
%    Asymptotic limits for Burgers equation}
%
% \author{\large \bf Paulo R. Zingano
% \thanks{This work was supported by CNPq (PQ 301236/91, AVG 7453835/96-1)
%         and FAPERGS (APC 96/1696-6)}
%         }
% \date{}
% \maketitle
%
%                     PAGE HEADINGS  settings:
%
% \thispagestyle{empty}
% \pagestyle{myheadings}
% \markboth{P. Zingano}{ASYMPTOTIC LIMITS FOR BURGERS EQUATION}
% \setlength{\baselineskip}{7mm}
%
%
% -----------------------------------------------------------------------
%
%
\begin{center}
{\large \bf
Some asymptotic limits for solutions of Burgers equation} \\
\nl
{\bf Paulo R. Zingano}$\mbox{}^{\dag}$ \\
Instituto de Matem\'atica \\
Universidade Federal do Rio Grande do Sul \\
Porto Alegre, RS 91509-900, Brazil
\nl
\end{center}
 \np {\bf \S 1 -- Introduction.}
In this paper,
we compute the limits \\
\nl
(1) \hspace{2.25cm}
$ {\displaystyle
{\gamma}_{\mbox{}_{\mbox{}_{\scriptstyle p}}}
\;=\; \limit \;\; \tpp \, \|\, u(\cdot,t) \,\nLp
, \hspace{0.50cm}
1 \leq p \leq \infty
} $ \\
\nl
\mbox{} \vspace{-0.35cm} \\
for solutions $u(\cdot,t)$ of the equation \\
\nl
(2) \hspace{4.0cm}
$ {\displaystyle
\ut \;+\; a \, \ux \;+\; b \, u \, \ux \;=\; c \, \uxx
} $ \\
\nl
satisfying the Cauchy condition \\
\nl
(3) \hspace{3.75cm}
$ {\displaystyle
u(x,0) \;=\; u_0(x), \hspace{0.50cm}
u_0 \in \Lone,
} $ \\
\nl
that is, $ \|\, u(\cdot,t) \,-\, u_0 \,\nLone \rightarrow\, 0\, $
as $ \, t \,\rightarrow\, 0 $, $ t > 0 $.
Here, $ \, \|\, u(\cdot,t) \,\nLp $ denotes the
\linebreak
$ L^p$ norm of $ \,u(\cdot,t)\, $
as a function of $ \,x\, $ for fixed $\,t$,
i.e., \\
\mbox{} \vspace{-0.20cm} \\
(4) \hspace{3.20cm}
$ {\displaystyle
\|\, u(\cdot,t) \,\nLp \;=\;
\left(\, \intline |\, u(x,t) \,|^{\mbox{}^{\scriptstyle p}} \; dx
\,\right)^{1/p}
} $ \\
\nl
if $ 1 \leq p < \infty $,
and \\
\nl
(5) \hspace{3.70cm}
$ {\displaystyle
\|\, u(\cdot,t) \,\nLinfty \;=\;\;\; \sup_{x \,\in\,{\bf R}}
\; |\, u(x,t) \,|
} $ \\
\nl
for $ p = \infty $.
In equation (2) above, $ a, b, c $ are real constants,
with $ c > 0 $. When $ b = 0 $ we have the familiar
heat equation; our main concern is the
case $ b \neq 0 $,
the so-called Burgers equation \mbox{[$\,1\,$]}, \mbox{[$\,3\,$]}.
Using the Hopf-Cole transformation \mbox{[$\,4\,$]}, \mbox{[$\,5\,$]},
it is well known that the solution in this case
is given by \\
\mbox{} \vspace{-0.30cm} \\
(6) \hspace{1.5cm}
$ {\displaystyle
u(x,t) \;=\; \factort \, \frac{1}{\,\varphi(x,t)\,} \,
\intline \expxyat
\, \phizero(y) \, u_0(y) \; dy
} $, \\
\nl
with \\
\mbox{} \vspace{-0.50cm} \\
(7) \hspace{2.0cm}
$ {\displaystyle
\varphi(x,t) \;=\; \factort \,
\intline \expxyat
\, \phizero(y) \; dy
} $, \\
\nl
\mbox{} \vspace{-0.3cm} \\
%
% ----------------------------------------------------------------------
\rule[-0.1cm]{10.35cm}{0.1mm} \\
\mbox{\footnotesize \mbox{}$^{\dag}\!$
\begin{minipage}[t]{10.0cm}
{\scriptsize
%This work was partly supported by CNPq (PQ 301236/91, AVG 7453835/96-1)
This work was partly supported by CNPq and Fapergs, Brazil. \\
E-mail address: pzingano@mat.ufrgs.br
}
\end{minipage}
}
% ----------------------------------------------------------------------
%
\newpage
\mbox{} \vspace{-0.6cm} \\
%
% ----------------------------------------------------------------------
%
where $ \, \phizero \in \Linfty \,$ is the Hopf-Cole transform of
the initial state $ \,u_0 $, i.e., \\
\mbox{} \vspace{-0.20cm} \\
(8) \hspace{4.0cm}
$ {\displaystyle
\phizero(x) \;=\; \HCxxizero \!
} $. \\
\nl
\mbox{} \vspace{-0.45cm} \\
From these expressions, we easily get that,
for every $ 1 \leq p \leq \infty $, one has \\
\mbox{} \vspace{-0.20cm} \\
(9) \hspace{4.0cm}
$ {\displaystyle
\|\, u(\cdot,t) \,\nLp \;\leq\; C \: \tmp
} $ \\
\nl
for some constant $ C > 0 $ which depends on the
magnitude of $ \,\|\, u_0 \,\|_{\mbox{}_{\scriptstyle 1}} $,
but it is not immediately evident how to compute
the limits $ {\gamma}_{\mbox{}_{\scriptstyle p}} $ above.
Denoting by $\,m\,$ the total mass of the solution, i.e., \\
\mbox{} \vspace{-0.20cm} \\
(10) \hspace{5.0cm}
$ {\displaystyle
m \;=\; \intline \!\!\! u_0(x) \; dx
} $, \\
\nl
we will show that \\
\mbox{} \vspace{-0.30cm} \\
(11) \hspace{2.0cm}
$ {\displaystyle
{\gamma}_{\mbox{}_{\scriptstyle p}} \;=\;
\frac{\;|\,m\,|}{\,\sqrt{\,4\,\pi\,c\,}\,} \;
(\, 4\,c\,)^{\mbox{}^{\scriptstyle
   \frac{\scriptstyle 1}{\scriptstyle 2\,p} }} \:
  \mbox{$ {\displaystyle \left| \frac{\mbox{}}{\mbox{}} \right. \! }$}
   \frac{\,2\,c\,}{\,b\,m\,} \, \mbox{\large $($} \, 1 \,-\, \ebm
   \! \mbox{\large $)$} \!
   \mbox{$ {\displaystyle \left. \frac{\mbox{}}{\mbox{}} \right| }$}
\: \|\, {\cal F} \,\nLp
} $ \\
\nl
\mbox{} \vspace{-0.30cm} \\
with $ \,{\cal F} \in \Lone \,\cap\, \Linfty \,$
being a function which depends on the parameters $ \,b, \,c, \,m $
above, given by \\
\mbox{} \vspace{-0.40cm} \\
(12) \hspace{4.0cm}
$ {\displaystyle
{\cal F} (x) \;=\; \frac{ \;\;\mbox{\sf e}^{\mbox{}^{\scriptstyle - \,
                          \mbox{\small $ x^2 $} }} }
                        { \, \mu \,-\, h \, \mbox{\sf erf}\,(x) \,}
} $, \\
\nl
\mbox{} \vspace{-0.40cm} \\
where $ \,\mbox{\sf erf}\,(x) \,$ is the error function \\
\nl
(13) \hspace{3.70cm}
$ {\displaystyle
\mbox{\sf erf}\,(x) \;=\;
\frac{1}{\,\sqrt{\,\pi\,}\;} \, \int_{0}^{x}
 \mbox{\sf e}^{\mbox{}^{\scriptstyle - \,
     \mbox{\footnotesize $\xi^{\mbox{}^{\scriptstyle 2}}$} }} \, d\xi
} $ \\
\nl
\mbox{} \vspace{-0.45cm} \\
and $ \mu $, $ h $
are given by \\
\mbox{} \vspace{-0.55cm} \\
(14) \hspace{2.20cm}
$ {\displaystyle
\mu \;=\; \frac{\;\; 1 \,+\, \ebm \!}{2}
} $,
\hspace{0.70cm}
$ {\displaystyle
h \;=\; |\;\, 1 \,-\, \ebm \!|
} $. \\
\nl
\mbox{} \vspace{-0.45cm} \\
When $ \, p = 1 $, (1), (11) become \\
\nl
(15) \hspace{4.0cm}
$ {\displaystyle
\limit \;\; \|\, u(\cdot,t) \,\nLone \;=\;\;\;
|\, m \,|
} $. \\
\nl
\mbox{} \vspace{-0.3cm} \\
In the case of heat equation, i.e., $ b = 0 $, the corresponding
results for $ {\gamma}_{\mbox{}_{\scriptstyle p}} $
are given by the limiting values of
the right-hand-side in (11) as $ \, b \,\tend\, 0 $, i.e., \\
\mbox{} \vspace{-0.20cm} \\
(16) \hspace{2.0cm} $ {\displaystyle
\limit \;\; \tpp \, \|\, u(\cdot,t) \,\nLp \;=\;
\frac{\;|\,m\,|}{\,\sqrt{\,4\,\pi\,c\,}\,} \,
\mbox{\Large $($}\, \frac{\,4\,\pi\,c\,}{p} \, \mbox{\Large $)$}
^{\mbox{}^{\mbox{$ \! \frac{\scriptstyle 1}{\scriptstyle \,2\,p\,} $}}}
                    } $ \\
\nl
\mbox{} \vspace{-0.30cm} \\
for every $ 1 \leq p \leq \infty $.
This case is easier and is briefly considered in Section 2 below.
The more interesting case when $ b \neq 0 $ is then taken up
in more detail in Section~3.
It is also shown in Section 3 that,
as it should be expected, the equations in the class (2)
are not asymptotically equivalent to one another:
if $ \, u, \hat{u} \,$ are solutions of \linebreak
\mbox{} \vspace{-0.15cm} \\
(17) \hspace{4.0cm}
$ {\displaystyle
\ut \;+\; a \, \ux \;+\; b \, u \, \ux \;=\; c \, \uxx
} $ \\
\mbox{} \vspace{-0.3cm} \\
and \\
\mbox{} \vspace{-0.30cm} \\
(18) \hspace{4.0cm}
$ {\displaystyle
\hat{u}_{\mbox{}_{\scriptstyle t}} \;+\;
\hat{a} \, \hat{u}_{\mbox{}_{\scriptstyle x}} \;+\;
\hat{b} \, \hat{u} \, \hat{u}_{\mbox{}_{\scriptstyle x}} \;=\;
\hat{c} \, \hat{u}_{\mbox{}_{\scriptstyle x x}}
} $ \\
\mbox{} \vspace{-0.1cm} \\
corresponding to the same initial profile $ \, u_0 \in \Lone \,$
with some nonzero mass,
and \linebreak
$ (\, a, \,b, \,c \,) \neq (\, \hat{a}, \,\hat{b}, \,\hat{c} \,) $,
then, for every $ 1 \leq p \leq \infty $,
there exist positive constants $ c_p , T_p $ such that \\
\mbox{} \vspace{-0.3cm} \\
(19) \hspace{3.50cm}
$ {\displaystyle
\|\, u(\cdot,t) \,-\, \hat{u}(\cdot,t) \,\nLp
\;\geq\;\; c_p \: \tmp
} $ \\
\nl
for all $ t \geq T_p $, so that
$ \, \|\, u(\cdot,t) - \hat{u}(\cdot,t) \,\|_{\mbox{}_{\scriptstyle p}} $
decays at exactly the same speed as
each term
$ \, \|\, u(\cdot,t) \,\|_{\mbox{}_{\scriptstyle p}} $,
$ \, \|\, \hat{u}(\cdot,t) \,\|_{\mbox{}_{\scriptstyle p}} $
on its own. \\
\nl
%
% ----------------------------------------------------------------------
%
%                             SECTION  2
%
%                      Heat  equation  (i.e., b = 0)
%
% -----------------------------------------------------------------------
%
\mbox{} \vspace{-0.3cm} \\
\np
{\bf \S 2 -- The case \mbox{\boldmath $ b = $} 0.}
Before we derive the results for the Burgers equation,
it will be convenient to consider briefly the simple case
of heat equation.
Clearly, it is sufficient to examine the case when $ a = 0 $,
so that we assume in this section that
$ \, u(\cdot,t) \in \Lone \cap \Linfty \,$
is the solution of the initial value problem \\
\nl
(20) \hspace{5.45cm}
$ {\displaystyle
\ut \;=\; c \, \uxx
} $ \\
\nl
(21) \hspace{4.95cm}
$ {\displaystyle
u(x,0) \;=\; u_0(x)
} $ \\
\nl
where $ \,c > 0\, $ is constant and $ u_0 \in \Lone $.
It is well known that $ u(x,t) $ is given by \\
\mbox{} \vspace{-0.30cm} \\
(22) \hspace{2.5cm}
$ {\displaystyle
u(x,t) \;=\;
\factort \intR \expxy  \, u_0(y) \; dy
} $, \\
\nl
\mbox{} \vspace{-0.30cm} \\
so that it satisfies,
for every $ 1 \leq p \leq \infty $, \\
\mbox{} \vspace{-0.10cm} \\
(23) \hspace{3.75cm}
$ {\displaystyle
\|\, u(\cdot,t) \,\nLp \;=\;
O(1) \: \tmp \!
} $. \\
\nl
A more subtle result which will be important throughout the
analysis is given in the following lemma. \\
\nl
{\bf Lemma 1} \\
{\em
Let $ \, u_0 \in \Lone \,$ be such that
$ {\displaystyle \intR \!\!\! u_0(x) \: dx \,=\, 0 }$.
Then, for every $ 1 \leq p \leq \infty $, one \linebreak
\mbox{} \vspace{-0.40cm} \\
has
} \\
\mbox{} \vspace{-0.30cm} \\
(24) \hspace{3.5cm}
$ {\displaystyle
\limit \;\; \tmp \: \|\, u(\cdot,t) \nLp \,=\; 0
} $. \\
\nl
\mbox{} \vspace{-0.70cm} \\
$\Box$
We will first show that
$ {\displaystyle
\, \limit \;\; \|\,u(\cdot,t)\,\nLone \;=\; 0
} $.
This has been shown for linear equations more general than (20) in
\mbox{[$\,2\,$]}, \mbox{[$\,6\,$]},
but for convenience we will give a \linebreak
\mbox{} \vspace{-0.4cm} \\
direct derivation below.
Given $ \,\varepsilon > 0 $, let $ A > 0 $ be chosen
such that \\
\mbox{} \vspace{-0.05cm} \\
\mbox{} \hspace{5.2cm}
$ {\displaystyle
\int_{\mbox{}_{\mbox{$ \scriptstyle |\,\mbox{\footnotesize $y$}\,|
                       \,\geq\,\mbox{\footnotesize $A$} $}}}
\!\!\!\!\!\! |\, u_0(y) \,| \; dy
\;\leq\; \varepsilon
} $, \\
\mbox{} \vspace{-0.05cm} \\
so that, from (22), \\
\mbox{} \vspace{-0.30cm} \\
\mbox{} \hspace{2.0cm}
$ {\displaystyle
\|\,u(\cdot,t)\,\nLone
\,\leq\;
\;\varepsilon \;+\,
\intline \!\! \factort \; \mbox{\Large $|$} \,
\int_{\mbox{}_{\mbox{$ \scriptstyle |\,\mbox{\footnotesize $y$}\,|
                       \,\leq\,\mbox{\footnotesize $A$} $}}}
\!\!\!\!\!\! \expxy  \, u_0(y) \; dy \;
\mbox{\Large $|$}  \; dx
} $ \\
\nl
\mbox{} \vspace{-0.30cm} \\
\mbox{} \hspace{3.75cm}
$ {\displaystyle
=\;
\;\varepsilon \;+\; \factorpi \!
\intline \mbox{\Large $|$} \, \int_{\mbox{}_{\mbox{$ \scriptstyle
                         |\,\mbox{\footnotesize $y$}\,|\,\leq\,
                             \mbox{\footnotesize $A$} $}}}
\!\!\!\!\!\! \expxiy  \!\! u_0(y) \; dy \; \mbox{\Large $|$} \; d\xi
} $. \\
\nl
\mbox{} \vspace{-0.10cm} \\
Letting $ \, t \,\tend\, \infty$, we then get \\
\mbox{} \vspace{-0.15cm} \\
\mbox{} \hspace{1.50cm}
$ {\displaystyle
\limitsup \;\; \|\,u(\cdot,t)\,\nLone
\;\leq\;
\;\varepsilon \;+\;
\factorpi \! \intline \! \expxi \! \mbox{\Large $|$} \,
\int_{\mbox{}_{\mbox{$ \scriptstyle |\,\mbox{\footnotesize $y$}\,|
                         \,\leq\, \mbox{\footnotesize $A$} $}}}
\!\!\!\!\!\!\!\! u_0(y) \; dy \; \mbox{\Large $|$} \; d\xi
} $ \\
\nl
\mbox{} \vspace{-0.30cm} \\
\mbox{} \hspace{4.0cm}
$ {\displaystyle
\leq\;
\;\varepsilon \;+\;
\mbox{\large $|$} \,
\int_{\mbox{}_{\mbox{$ \scriptstyle |\,\mbox{\footnotesize $y$}\,|
                                \,\leq\, \mbox{\footnotesize $A$} $}}}
\!\!\!\!\!\!\!\! u_0(y) \; dy \;
\mbox{\Large $|$}
\;\leq \;
2 \, \varepsilon
} $, \\
\nl
\mbox{} \vspace{-0.30cm} \\
where we have used that
$ {\displaystyle \intline \!\!\!\! u_0(y) \, dy \,=\, 0 }$.
Since $ \,\varepsilon > 0\, $ is arbitrary, this gives
\[
\limit \;\; \|\,u(\cdot,t) \,\nLone \;=\; 0,
\]
which concludes the case $ p = 1 $.
Now, given $ \, 1 < p < \infty $, we have, using (23), \\
\mbox{} \vspace{-0.20cm} \\
\mbox{} \hspace{0.50cm}
$ {\displaystyle
\tpp \, \|\, u(\cdot,t) \,\nLp
\;\leq\;
\; \|\, u(\cdot,t) \,\nLoneINTP \,
\mbox{\Large $($}\,
\thalf \, \|\, u(\cdot,t) \,\nLinfty
\mbox{\Large $)$}^{\mbox{}^{\scriptstyle 1 \,-\,
   \frac{\scriptstyle 1}{\scriptstyle p} }}
\leq\;
C \: \|\, u(\cdot,t) \,\nLoneINTP
} $ \\
\nl
for some constant $ \, C > 0 $, so that
\[
\limit \;\; \tpp \, \|\, u(\cdot,t) \,\nLp \;=\; 0
\]
from the previous case.
Finally, we consider the case $ \, p = \infty $:
from (22), it readily follows that
\[
\|\, \ux(\cdot,t) \,\nLtwo
\;=\;
O(1) \: t^{\mbox{}^{\scriptstyle - \frac{\scriptstyle 3}{\scriptstyle 4} }},
\]
so that
\[
\thalf \, \|\, u(\cdot,t) \,\nLinfty
=\;
O(1) \: \thalf \, \|\, u(\cdot,t) \,\nLtwohalf \:
\|\, u_{\mbox{}_{\scriptstyle x}}(\cdot,t) \,\nLtwohalf
\;=\;
O(1) \, \mbox{\Large $($}\, \thalf \, \|\, u(\cdot,t) \,\nLtwo
\,\mbox{\Large $)$}^{\mbox{}^{\scriptstyle
                        \frac{\scriptstyle 1}{\scriptstyle 2} }}
\!,
\]
\nl
which gives the result from the case $ \, p = 2\, $ already considered.
\hfill $\Box$ \\
%
% ---------------------------------------  END OF PROOF  for Lemma 1
%
\np
Using the estimates above,
we can easily obtain the limits $ \,{\gamma}_{\mbox{}_{\scriptstyle p}} \, $
for (20), (21), as shown next. \\
\nl
{\bf Theorem 1} \\
{\em
Given $ \,u_0 \in \Lone $, the solution $ \,u(x,t) $ of\/ $(20)$, $(21)$
satisfies \\
\nl
$(i)$
\hspace{0.40cm}
$ {\displaystyle
\limit \;\; \|\, u(\cdot,t) \,\nLone \;=\;\; |\,m\,|
} $, \\
\nl
$(ii)$
\hspace{0.30cm}
$ {\displaystyle
\limit \;\; \tpp \, \|\, u(\cdot,t) \,\nLp \;=\;\;
\frac{\;|\,m\,|}{\,\sqrt{\,4\,\pi\,c\,}\,} \, \mbox{\Large $($}\,
\frac{ \,4\,\pi\,c\, }{ p } \,
\mbox{\Large $)$}^{\mbox{}^{\scriptstyle
          \frac{\scriptstyle 1}{\scriptstyle 2\,p} }}
} $, \\
\mbox{} \vspace{-0.30cm} \\
and \\
\mbox{} \vspace{-0.30cm} \\
$(iii)$
\hspace{0.20cm}
$ {\displaystyle
\limit \;\; \thalf \, \|\, u(\cdot,t) \,\nLinfty \;=\;
\frac{ |\,m\,| }{\, \sqrt{\,4\,\pi\,c\,}\,}
} $, \\
\nl
\mbox{} \vspace{-0.35cm} \\
where $ {\displaystyle \; m \,=\, \intR \!\!\!\! u_0(x) \;dx }$.
} \\
\nl
\mbox{} \vspace{-0.15cm} \\
$\Box$ In fact, this can be readily established for \\
\mbox{} \vspace{-0.30cm} \\
\mbox{} \hspace{4.0cm}
$ {\displaystyle
\tilde{u}(x,t) \;=\; \factortm \, \int_{0}^{1} \expxy \; dy
} $, \\
\nl
\mbox{} \vspace{-0.30cm} \\
i.e., the solution of $\,$(20), (21) $\,$corresponding
to the initial profile
$ \,\tilde{u}_0(x) = m \,
\mbox{\tiny $\mbox{}^{\mbox{\normalsize $\chi$}}$}
   _{\mbox{}_{\scriptstyle [\,0,\,1\,]}}(x) $,
where
$ \mbox{\tiny $\mbox{}^{\mbox{\normalsize $\chi$}}$}
   _{\mbox{}_{\scriptstyle [\,0,\,1\,]}} $
denotes the
characteristic function of the interval $ \, [\,0, 1 \,] $.
The result then follows for an arbitrary $ \, u_0 \in \Lone $ with
the same mass $\,m$, because, from Lemma~1, one has \\
\mbox{} \vspace{-0.30cm} \\
\mbox{} \hspace{3.50cm}
$ {\displaystyle
\limit \;\; \tpp \, \|\, u(\cdot,t) \,-\, \tilde{u}(\cdot,t) \,\nLp \;=\; 0
} $. \hfill $\Box$ \\
\nl
%
% ----------------------------------------------------------------------
%
%                             SECTION  3
%
%                          Burgers  equation
%
% -----------------------------------------------------------------------
%
\mbox{} \vspace{-0.30cm} \\
\np
{\bf \S 3 -- The case \mbox{\boldmath $ b \neq $} 0.}
In this section we extend the analysis above
to the more interesting case of Burgers equation.
As before, we assume without loss of generality $\, a = 0 $,
and let $ \, u(\cdot,t) \,$ be the solution
of the Cauchy problem \\
\nl
(25a) \hspace{4.0cm}
$ {\displaystyle
\ut \;+\; b \, u \, \ux \;=\; c \, \uxx
} $, \\
\nl
(25b) \hspace{4.2cm}
$ {\displaystyle
u(x,0) \;=\; u_0(x)
} $, \\
\nl
where $\, b \neq 0 $, $ c > 0 $, and $ \, u_0 \in \Lone $.
Using the Hopf-Cole transformation \mbox[$\,3\,$], \mbox[$\,4\,$] \\
\mbox{} \vspace{-0.20cm} \\
(26) \hspace{4.0cm}
$ {\displaystyle
\varphi(x,t) \;=\; \HCxxit \!
} $, \\
\nl
\mbox{} \vspace{-0.35cm} \\
we have that $ \,\varphi(\cdot,t) \,$ satisfies \\
(27a) \hspace{4.75cm}
$ {\displaystyle
\phit \;=\; c \, \phixx
} $ \\
\mbox{} \vspace{-0.40cm} \\
(27b) \hspace{3.75cm}
$ {\displaystyle
\varphi(x,0) \;=\; \phizero(x) \;\equiv\; \HCxxizero
} $ \\
\nl
with $ \,u(x,t)\, $ given by \\
\mbox{} \vspace{-0.20cm} \\
(28) \hspace{3.75cm}
$ {\displaystyle
u(x,t) \;=\; - \, \frac{\,2\,c\,}{b} \:
\frac{\,\phix(x,t)\,}{\varphi(x,t)}
} $, \\
\nl
that is, \\
\mbox{} \vspace{-0.50cm} \\
(29) \hspace{1.50cm}
$ {\displaystyle
u(x,t) \;=\; \factort \, \frac{1}{\,\varphi(x,t)\,} \,
\intline \expxy \, \phizero(y) \, u_0(y) \; dy
} $. \\
\nl
\mbox{} \vspace{-0.10cm} \\
It follows from this expression that
$ \,u(\cdot,t) \in \Lone \cap \Linfty \,$
for all $ \,t > 0 $,
with \\
\nl
(30) \hspace{4.0cm}
$ {\displaystyle
\|\, u(\cdot,t) \,\nLp \,=\; O(1) \: \tmp
} $ \\
\nl
for every $\, 1 \leq p \leq \infty $.
Moreover, when $ \, u_0 \,$ has zero mass,
we readily get the following estimate from Lemma 1. \\
\nl
{\bf Lemma 2} \\
{\em
Let $ \, u_0 \in \Lone \, $ be such that
$ {\displaystyle \intR \!\!\! u_0(x) \, dx \,=\, 0 }$.
Then, for every $ 1 \leq p \leq \infty $, one \linebreak
\mbox{} \vspace{-0.50cm} \\
has
} \\
\mbox{} \vspace{-0.40cm} \\
(31) \hspace{3.75cm}
$ {\displaystyle
\limit \;\; \tmp \, \|\, u(\cdot,t) \nLp \,=\; 0
} $. \\
\nl
\mbox{} \vspace{-0.30cm} \\
$\Box$
In fact, from (27) we see that $\,\phix\,$ satisfies
the conditions of Lemma 1, so that, for every $ \,1 \leq p \leq \infty $,\\
\mbox{} \vspace{-0.30cm} \\
\mbox{} \hspace{4.00cm}
$ {\displaystyle
\limit \;\; \tpp \: \|\, \phix(\cdot,t) \,\nLp \;=\; 0
} $. \\
\nl
Since $\, 1/\varphi(\cdot,t) \,$ is uniformly bounded,
the same is true of
$\, u(\cdot,t) \,$ in view of (28).
\hfill $\Box$ \\
\nl
Another fundamental consequence of Lemma 1 is given next. \\
\nl
% -------------------------------------------------------- Lemma 3
{\bf Lemma 3} \\
{\em
Let $ \,u_0, v_0 \in \Lone \,$ be such that
} \\
\mbox{} \vspace{-0.20cm} \\
(32) \hspace{3.5cm}
$ {\displaystyle
\intline \!\!\! u_0(x) \; dx \;=\; \intline \!\!\! v_0(x) \; dx
} $, \\
\nl
\mbox{} \vspace{-0.4cm} \\
{\em
and let $ \, u(\cdot,t) $, $ v(\cdot,t) \,$ be the solutions of
$\,(25)$ corresponding to the initial profiles $ \,u_0 $, $ v_0 $,
respectively.
Then, for every $ 1 \leq p \leq \infty $, one has
} \\
\mbox{} \vspace{-0.20cm} \\
(33) \hspace{2.75cm}
$ {\displaystyle
\limit \;\; \tpp \, \|\, u(\cdot,t) \,-\, v(\cdot,t) \,\nLp \,=\; 0
} $. \\
\nl
% -------------------------------------------------------------------
%
\mbox{} \vspace{-0.15cm} \\
$\Box$
Letting $ \,\varphi(\cdot,t)$, $ \psi(\cdot,t) \,$
be the Hopf-Cole transforms of $ \,u(\cdot,t)$, $ v(\cdot,t)$,
respectively, i.e.,
\[
\varphi(x,t) \;=\; \HCxxit
\!, \hspace{0.70cm}
\psi(x,t) \;=\; \HCvxxit
\!,
\]
\mbox{} \vspace{-0.30cm} \\
and setting
$ \,\omega = \phix - \psix $,
we have that $ \omega(\cdot,t) $ has zero mass and satisfies
$\,\omegat \;=\; c \, \omegaxx$,
so that, from Lemma 1, \\
\mbox{} \vspace{-0.30cm} \\
\mbox{} \hspace{3.00cm}
$ {\displaystyle
\limit \;\; \tpp \: \|\, \phix(\cdot,t) \,-\, \psix(\cdot,t) \,\nLp \,=\; 0
} $ \\
\nl
for every $ \,1 \leq p \leq \infty $,
i.e., \\
\mbox{} \vspace{-0.20cm} \\
\mbox{} \hspace{2.00cm}
$ {\displaystyle
\limit \;\; \tpp \:
\|\, \varphi(\cdot,t) \, u(\cdot,t) \,-\, \psi(\cdot,t) \, v(\cdot,t) \,\nLp
\,=\; 0
} $. \\
\nl
Since $ \,1/\varphi(\cdot,t) $, $ 1/\psi(\cdot,t) \,$
are uniformly bounded
and \\
\mbox{} \vspace{-0.20cm} \\
\mbox{} \hspace{4.00cm}
$ {\displaystyle
\limit \;\; \|\, \varphi(\cdot,t) \,-\, \psi(\cdot,t) \,\nLinfty
\,=\; 0
} $, \\
\nl
we get the result.
\hfill $\Box$ \\
% --------------------------------------------- End of PROOF for Lemma 3
\nl
We are now in position to compute the limits
$\,{\gamma}_{\mbox{}_{\mbox{}_{\scriptstyle p}}}$
for an arbitrary $ \,u_0 \,$ in $ \,\Lone $. \\
\nl
%
% ------------------------------------------------------------ Theorem 2
{\bf Theorem 2} \\
{\em
Given $ \,u_0 \in \Lone $, the solution $ \,u(\cdot,t) $ of\/ $(25)$
satisfies \\
\nl
$(i)$
\hspace{0.40cm}
$ {\displaystyle
\limit \;\; \|\, u(\cdot,t) \,\nLone \;=\;\; |\,m\,|
} $, \\
\nl
$(ii)$
\hspace{0.30cm}
$ {\displaystyle
\limit \;\; \tpp \, \|\, u(\cdot,t) \,\nLp \;=\;\;
\frac{\;|\,m\,|}{\,\sqrt{\,4\,\pi\,c\,}\,} \,
(\, 4\,c\,)^{\mbox{}^{\scriptstyle
              \frac{\scriptstyle 1}{\scriptstyle 2\,p} }} \,
\left|\, \frac{\,2\,c\,}{\,b\,m\,} \, \mbox{\Large $($}\, 1 \,-\, \ebm
\,\mbox{\Large $)$} \,\right| \:
\|\, {\cal F} \,\nLp
} $ \\
\mbox{} \vspace{-0.35cm} \\
and \\
$(iii)$
\hspace{0.20cm}
$ {\displaystyle
\limit \;\; \thalf \, \|\, u(\cdot,t) \,\nLinfty \;=\;
\frac{\;|\,m\,|}{\,\sqrt{\,4\,\pi\,c\,}\,} \,
\left|\, \frac{\,2\,c\,}{\,b\,m\,} \, \mbox{\Large $($} \, 1 \,-\, \ebm
\,\mbox{\Large $)$} \,\right| \:
\|\, {\cal F} \,\nLinfty
} $, \\
\nl
\mbox{} \vspace{-0.30cm} \\
where $ \,{\cal F} \in \Lone \cap \Linfty \,$ is given in
$\,(12) - (15)$,
and
$ {\displaystyle \:m \,=\, \intR \!\!\! u_0(x) \;dx }$.
} \\
\nl
% ---------------------------------------------------------- Theorem 2
\mbox{} \vspace{-0.3cm} \\
$\Box$
Because of Lemma 3, it is sufficient to show the result for the
particular initial state
$\, u_0 = m\, \mbox{\tiny $\mbox{}^{\mbox{\normalsize $\chi$}}$}
_{[\,0,\,1\,]} $,
in which case $\, u(\cdot,t) \, $ is given by \\
\mbox{} \vspace{-0.30cm} \\
(34) \hspace{2.0cm}
$ {\displaystyle
u(x,t) \;=\; \factormt \, \frac{1}{\,\varphi(x,t)\,} \,
\int_{0}^{1}
\expxyat \, \phizero(y) \; dy
} $, \\
\nl
\mbox{} \vspace{-0.35cm} \\
where \\
\mbox{} \vspace{-0.35cm} \\
(35) \hspace{2.50cm}
$ {\displaystyle
\varphi(x,t) \;=\; \factort \,
\intline
\expxyat \, \phizero(y) \; dy
} $, \\
\nl
\mbox{} \vspace{-0.40cm} \\
with \\
\mbox{} \vspace{-0.50cm} \\
(36) \hspace{4.00cm}
$ {\displaystyle
\phizero(x) \;=\; \HCxxizero
} $, \\
\nl
\mbox{} \vspace{-0.35cm} \\
see (6), (7), (8).
In particular, for any $\, t > 0 $, \\
\mbox{} \vspace{-0.10cm} \\
\mbox{} \hspace{2.00cm}
$ {\displaystyle
\|\, u(\cdot,t) \,\nLone
\,=\;
\left|\; \intline \!\!\! u(x,t) \; dx \;\right|
\;=\;
\left|\; \intline \!\!\!u_0(x) \; dx \;\right|
\;=\;
|\, m \,|
} $, \\
\nl
\mbox{} \vspace{-0.30cm} \\
which shows $(i)$.
To get $(ii)$, $(iii)$, we introduce \\
\mbox{} \vspace{-0.30cm} \\
(37) \hspace{3.50cm}
$ {\displaystyle
\Hzerox \;=\; \left\{\, \begin{array}{lll}
                          1\,, & \mbox{} & x < \alpha \\
                          \mbox{} \vspace{-0.30cm} \\
                          \ebm \!\!, & \mbox{} & x > \alpha
                        \end{array}
              \right.
} $ \\
\nl
where $\, \alpha > 0 \,$ is chosen so that \\
\nl
(38) \hspace{3.50cm}
$ {\displaystyle
\intline (\, \Hzerox \,-\, \phizero(x) \,) \; dx \;=\; 0
} $, \\
\nl
i.e., \\
\mbox{} \vspace{-0.40cm} \\
(39) \hspace{2.00cm}
$ {\displaystyle
\alpha \;+\; (\, 1 \,-\, \alpha \,) \: \ebm
\;=\;
\frac{\,2\,c\,}{\,b\,m\,} \: (\, 1 \,-\, \ebm \!)
} $, \\
\nl
as illustrated in the picture below. \\
%
% -----------------------------------------------------------
%                            FIGURA 1
% -----------------------------------------------------------
%
\begin{center}
\setlength{\unitlength}{1.0cm}
\begin{picture}(10.000, 2.500)
%
% ----------------------------------------------
%
%      ( x , u )   --->   ( xPLOT , uPLOT )
%
%                     x  -  xMIN
%        xPLOT  =  ---------------  Width
%                    xMAX - xMIN
%
%                     u  -  uMIN
%        uPLOT  =  ---------------  Height
%                    uMAX - uMIN
%
%        xMIN  =    -.250000E+01     xMAX  =     .275000E+01
%        uMIN  =    -.750000E-01     uMAX  =     .130000E+01
%
%        Width =     .100000E+02    Height =     .250000E+01
%
% ----------------------------------------------
%
%  Specifying the coordinate FRAME:
%
%    ( reference axes )
%
%     xINITIAL = -.200000E+01  xFINAL =  .200000E+01
%     number of subintervals in x-direction =   4
%
%     uINITIAL =  .000000E+00  uFINAL =  .100000E+01
%     number of subintervals in u-direction =   1
%
%     Location of coordinate axes:
%        x-axis:  u =  .000000E+00
%        u-axis:  x =  .000000E+00
%
\put (  0.000 ,  0.146 )   { \line(1,0) { 10.100 } }
\put (  0.000 ,  0.136 ) { \vector(1,0) { 10.200 } }
\put (  0.000 ,  0.126 )   { \line(1,0) { 10.100 } }
\put (  4.752 ,  0.000 )   { \line(0,1) {  2.585 } }
\put (  4.762 ,  0.000 ) { \vector(0,1) {  2.675 } }
\put (  4.772 ,  0.000 )   { \line(0,1) {  2.585 } }
\put ( 10.300 ,  0.136 ) { \makebox(0,0)[l] {$x$} }
%\put ( 4.762 ,  2.775 ) { \makebox(0,0)[b] {$u$} }
%
\put (  0.952 ,  0.036 )  { \line(0, 1) {0.20} }
\put (  0.952 , -0.100 ) { \makebox(0,0)[t] { -2 } }
\put (  2.857 ,  0.036 )  { \line(0, 1) {0.20} }
\put (  2.857 , -0.100 ) { \makebox(0,0)[t] { -1 } }
\put (  4.762 ,  0.036 )  { \line(0, 1) {0.20} }
\put (  4.762 , -0.100 ) { \makebox(0,0)[t] {  0 } }
\put (  6.667 ,  0.036 )  { \line(0, 1) {0.20} }
\put (  6.667 , -0.100 ) { \makebox(0,0)[t] {  1 } }
\put (  8.571 ,  0.036 )  { \line(0, 1) {0.20} }
\put (  8.571 , -0.100 ) { \makebox(0,0)[t] {  2 } }
\put (  4.662 ,  0.136 )  { \line( 1,0) {0.20} }
%\put (  4.562 ,  0.136 ) { \makebox(0,0)[r] { 0 } }
\put (  4.662 ,  1.955 )  { \line( 1,0) {0.20} }
%\put (  4.562 ,  1.955 ) { \makebox(0,0)[r] { 1 } }
%
% -----------------------------------------------------------------------
%
%                  Function  Hzero(x):
%
% -----------------------------------------------------
%
\put ( 5.5581 , 1.9545 ) { \line(-1,0) {5.0819} }
\put ( 5.5581 , 0.8052 ) { \line( 1,0) {3.9657} }
\put ( 5.5581 , 0.8052 ) { \line( 0,1) {1.1493} }
\put ( 5.5581 , 0.0360 ) { \line( 0,1) {0.20} }
\multiput ( 5.5581 , 0.285 )( 0 , 0.173 ) { 3 }{ \line(0,1) {0.087} }
\multiput ( 6.6667 , 0.285 )( 0 , 0.173 ) { 3 }{ \line(0,1) {0.087} }
\put ( 4.662 , 0.805 ) { \line( 1,0) {0.20} }
\multiput ( 4.962 , 0.805 )( 0.200 , 0 ) { 3 }{ \line(1,0) {0.100} }
\put ( 4.7619 , 1.9545 ) { \circle*{0.07} }
\put ( 6.6667 , 0.8052 ) { \circle*{0.07} }
\put ( 4.7600 , 2.0500 ) { \makebox(0,0)[br] { 1 } }
\put ( 4.9000 , 0.8052 ) { \makebox(0,0)[br] {\mbox{$ \ebm $}} }
\put ( 5.5581 , -0.150 ) { \makebox(0,0)[t] {\mbox{$ \alpha $}} }
\put ( 5.9000 , 1.3500 ) { \makebox(0,0)[bl]
       {\mbox{$ {\varphi}_{\mbox{}_{\mbox{}_{\scriptstyle \!0}}} $}} }
\put ( 2.5000 , 2.1000 ) { \makebox(0,0)[b] {\mbox{$ {\cal H}_0 $}} }
\put ( 8.0000 , 0.9200 ) { \makebox(0,0)[b] {\mbox{$ {\cal H}_0 $}} }
\put ( 4.7600 , -0.800 ) { \makebox(0,0)[t]
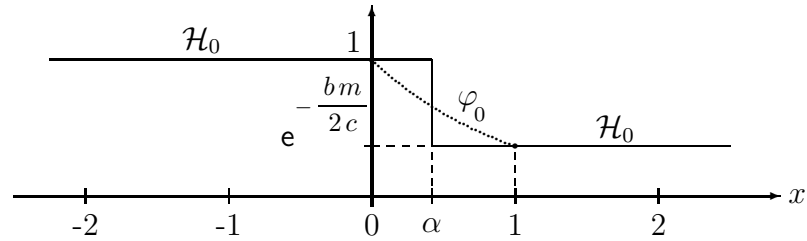
       {\mbox{\sf Figure 1: $ \,{\cal H}_{\mbox{}_{\scriptstyle 0}} \,$
                    and $ \,{\varphi}_{\mbox{}_{\scriptstyle 0}} $}} }
%
% ------------------------------------------------------------------------
%
%                      Function  Phizero(x):
%
% ------------------------------------------------
%
\put (   4.7619 ,   1.9545 ){ \circle*{0.05} }
\put (   4.8055 ,   1.9134 ){ \circle*{0.05} }
\put (   4.8496 ,   1.8727 ){ \circle*{0.05} }
\put (   4.8942 ,   1.8325 ){ \circle*{0.05} }
\put (   4.9393 ,   1.7929 ){ \circle*{0.05} }
\put (   4.9848 ,   1.7538 ){ \circle*{0.05} }
\put (   5.0307 ,   1.7152 ){ \circle*{0.05} }
\put (   5.0772 ,   1.6772 ){ \circle*{0.05} }
\put (   5.1240 ,   1.6397 ){ \circle*{0.05} }
\put (   5.1714 ,   1.6029 ){ \circle*{0.05} }
\put (   5.2191 ,   1.5665 ){ \circle*{0.05} }
\put (   5.2673 ,   1.5308 ){ \circle*{0.05} }
\put (   5.3159 ,   1.4957 ){ \circle*{0.05} }
\put (   5.3650 ,   1.4611 ){ \circle*{0.05} }
\put (   5.4145 ,   1.4271 ){ \circle*{0.05} }
\put (   5.4643 ,   1.3938 ){ \circle*{0.05} }
\put (   5.5146 ,   1.3610 ){ \circle*{0.05} }
\put (   5.5653 ,   1.3289 ){ \circle*{0.05} }
\put (   5.6163 ,   1.2974 ){ \circle*{0.05} }
\put (   5.6677 ,   1.2664 ){ \circle*{0.05} }
\put (   5.7195 ,   1.2361 ){ \circle*{0.05} }
\put (   5.7716 ,   1.2064 ){ \circle*{0.05} }
\put (   5.8241 ,   1.1774 ){ \circle*{0.05} }
\put (   5.8769 ,   1.1489 ){ \circle*{0.05} }
\put (   5.9301 ,   1.1210 ){ \circle*{0.05} }
\put (   5.9835 ,   1.0938 ){ \circle*{0.05} }
\put (   6.0373 ,   1.0671 ){ \circle*{0.05} }
\put (   6.0914 ,   1.0411 ){ \circle*{0.05} }
\put (   6.1457 ,   1.0156 ){ \circle*{0.05} }
\put (   6.2003 ,   0.9908 ){ \circle*{0.05} }
\put (   6.2552 ,   0.9665 ){ \circle*{0.05} }
\put (   6.3103 ,   0.9428 ){ \circle*{0.05} }
\put (   6.3657 ,   0.9197 ){ \circle*{0.05} }
\put (   6.4213 ,   0.8972 ){ \circle*{0.05} }
\put (   6.4771 ,   0.8752 ){ \circle*{0.05} }
\put (   6.5332 ,   0.8538 ){ \circle*{0.05} }
\put (   6.5894 ,   0.8329 ){ \circle*{0.05} }
\put (   6.6459 ,   0.8126 ){ \circle*{0.05} }
%
% ---------------------------------------------------
%
\end{picture}
\end{center}
%
% ----------------------------------------------------- END OF FIGURE 1
%
\np
\nl
\nl
\mbox{} \vspace{-0.35cm} \\
Setting \\
\mbox{} \vspace{-0.50cm} \\
(40) \hspace{2.50cm}
$ {\displaystyle
{\cal H}(x,t) \;=\;
\factort \, \intline \expxy \, \Hzeroy \; dy
} $, \\
\nl
\mbox{} \vspace{-0.30cm} \\
we have, from (38) and Lemma 1, \\
\mbox{} \vspace{-0.10cm} \\
(41) \hspace{2.75cm}
$ {\displaystyle
\limit \;\; \tpp \:
\|\, {\cal H}(\cdot,t) \,-\, \varphi(\cdot,t) \,\nLp
\;=\; 0
} $ \\
\nl
\mbox{} \vspace{-0.40cm} \\
for every $\, 1 \leq p \leq \infty $,
so that \\
\mbox{} \vspace{-0.10cm} \\
(42) \hspace{2.90cm}
$ {\displaystyle
\limit \;\; \tpp \:
\|\, u(\cdot,t) \,-\, \omega(\cdot,t) \,\nLp
\;=\; 0
} $, \\
\nl
\mbox{} \vspace{-0.40cm} \\
where $\,\omega(\cdot,t)\,$ is defined by \\
\mbox{} \vspace{-0.20cm} \\
(43) \hspace{1.50cm}
$ {\displaystyle
\omega(x,t) \;=\; \factormt \, \frac{1}{\,{\cal H}(x,t)\,} \,
\int_{0}^{1} \mbox{\sf e}^{\mbox{}^{\scriptstyle - \,
\frac{\mbox{$\scriptstyle \:(\, x - y - a\,t \,)^2$}}
     {\mbox{$\scriptstyle 4\,c\,t$}} }}\,
\Hzeroy \; dy
} $, \\
\nl
\mbox{} \vspace{-0.30cm} \\
and $ \,{\cal H}(\cdot,t) \,$ is given in (40) above,
that is, \\
\mbox{} \vspace{-0.0cm} \\
(44) \hspace{3.5cm}
$ {\displaystyle
{\cal H}(x,t) \;=\; \mu \;-\; \sigma \, h \: \mbox{\sf erf}
\,\mbox{\Large $($}\, \frac{\: x\,-\, \alpha \,}{\,\sqrt{\,4\,c\,t\,}\,}
\,\mbox{\large $)$}
} $, \\
\nl
\mbox{} \vspace{-0.35cm} \\
where $\,\sigma\,$ is the sign of the product $\, b\,m \,$
(i.e., $ \,\sigma = 1\, $ if $ \, b\,m > 0$, $ \sigma = -\,1\,$
otherwise) and $\,\mu$, $h$, $\mbox{\sf erf}\,(x)\,$ are given
in $(13) - (15)$.
We will now derive $(ii)$, for $ \, 1 \leq p < \infty $:
given $\,\xi \in {\bf R} $, we have \\
\mbox{} \vspace{-0.30cm} \\
\mbox{} \hspace{5.25cm}
$ {\displaystyle
\omega \,(\, \alpha \,+\, \xi \, \sqrt{\,4\,c\,t\,} \,, \,t\,) \;} $ \\
(45) \\
\mbox{} \vspace{-0.80cm} \\
\mbox{} \hspace{2.50cm}
$ {\displaystyle
=\;
\factormt \: \frac{1}{\,\mu \,-\, \sigma \, h \, \mbox{\sf erf}\,(\xi) \,}
\, \int_{0}^{1} \expxiya \, \Hzeroy \; dy
} $, \\
\nl
so that \\
\mbox{} \vspace{-0.30cm} \\
\mbox{} \hspace{5.25cm}
$ {\displaystyle
\tppp \: \|\, \omega(\cdot,t) \,\nLpp
\;} $ \\
\mbox{} \vspace{-0.55cm} \\
$ {\displaystyle
=\;
\mbox{\Large $($}\, \factormm \,\mbox{\Large $)$}^{\mbox{}^{\scriptstyle p}}
\sqrt{\,4\,c\;} \, \intline
\frac{1}{\;\left|\; \mu \,-\, \sigma \, h \, \mbox{\sf erf}\,(\xi) \;
         \right|^{\mbox{}^{\scriptstyle p}} }
\:
\mbox{\huge $|$} \, \int_{0}^{1} \expxiya \!\! \Hzeroy \; dy
\;\mbox{\huge $|$}^{\mbox{}^{\scriptstyle p}}
\; d\xi
} $. \\
\nl
\mbox{} \vspace{-0.0cm} \\
Since, for all $\, \xi \in {\bf R} \,$ and
$ \, t \geq 1/4c $, we have \\
\mbox{} \vspace{-0.10cm} \\
\mbox{} \hspace{2.0cm}
$ {\displaystyle
\mbox{\huge $|$} \, \int_{0}^{1} \expxiya \!\! \Hzeroy \; dy
\;\mbox{\huge $|$}^{\mbox{}^{\scriptstyle p}}
\;\leq\;
\;\mbox{\sf e}^{\mbox{}^{\scriptstyle -\,
\mbox{\footnotesize $ {\displaystyle
\frac{\,\mbox{\footnotesize $p$} \,}{\,\mbox{\footnotesize $2$} \,}
\, {\xi}^{\mbox{}^{\scriptstyle 2}} \,+\, 1 }$} }}
\, \|\, \Hzero \,\|^{\mbox{}^{\scriptstyle p}}
                   _{\mbox{}_{\mbox{}_{\scriptstyle L^{1}(0,1)}}}
} $, \\
\nl
\mbox{} \vspace{-0.30cm} \\
we get, by Lebesgue's dominated convergence theorem, \\
\mbox{} \vspace{-0.10cm} \\
\mbox{} \hspace{4.50cm}
$ {\displaystyle
\limit \;\; \tpp \, \|\, \omega(\cdot,t) \,\nLp \,} $ \\
\mbox{} \vspace{-0.05cm} \\
\mbox{} \hspace{0.50cm}
$ {\displaystyle
=\;
\factormm \, \fourcp
\left(\, \int_{0}^{1} \Hzeroy \; dy \,\right) \:
\left(\, \intline \mbox{\huge $|$} \,
\frac{\expxi}{\,\mu\,-\,\sigma\,h\,\mbox{\sf erf}\,(\xi)\,}
\,\mbox{\huge $|$}^{\mbox{}^{\scriptstyle p}}
\; d\xi \;\right)^{\!1/p}
} $ \\
\mbox{} \vspace{-0.05cm} \\
\mbox{} \hspace{2.50cm}
$ {\displaystyle
=\;
\factormm \, \fourcp \, \fbcm \, \|\, {\cal F} \,\nLp
} $ \\
\nl
\mbox{} \vspace{-0.20cm} \\
in view of (37), (39). This shows $(ii)$.
Finally, for $ \, p = \infty $, we observe that,
letting $ \, t \,\tend\,\infty \,$ in (45), we get \\
\mbox{} \vspace{-0.20cm} \\
\mbox{} \hspace{0.50cm}
$ {\displaystyle
\limitinf \;\; \thalf \: \|\, \omega(\cdot,t) \,\nLinfty
\;\geq\;
\factormm \: \fbcm \: \frac{ \, \mbox{\sf e}^{\mbox{}^{\scriptstyle
       -\, \mbox{\footnotesize $ {\xi}^{\mbox{}^{\scriptstyle 2}} $} }} }
                           {\, \mu \,-\, \sigma\,h\,\mbox{\sf erf}\,(\xi)\,}
} $ \\
\nl
\mbox{} \vspace{-0.30cm} \\
for every $ \,\xi \in {\bf R} $,
so that \\
\mbox{} \vspace{-0.40cm} \\
(46) \hspace{0.75cm}
$ {\displaystyle
\limitinf \;\; \thalf \: \|\, \omega(\cdot,t) \,\nLinfty
\;\geq\;
\factormm \: \fbcm \: \|\, {\cal F} \,\nLinfty
} $. \\
\nl
\mbox{} \vspace{-0.10cm} \\
On the other hand, for $ \, t > 0 \, $ let
$ \,{\xi}_{\mbox{}_{\mbox{}_{\scriptstyle t}}} \in {\bf R} \,$
be such that \\
\nl
(47) \hspace{3.25cm}
$ {\displaystyle
\|\: \omega(\cdot,t)\,\nLinfty \;=\;
|\; \omega(\, \alpha \,+\, {\xi}_{\mbox{}_{\mbox{}_{\scriptstyle t}}}
  \, \sqrt{\, 4\,c\,t\,} \,,\, t \,) \;|
} \,$; \\
\nl
since
$ \, \limitinf \;\, \thalf \, \|\: \omega(\cdot,t) \,\nLinfty > 0 \,$
from (46), we must have
$ \, {\xi}_{\mbox{}_{\mbox{}_{\scriptstyle t}}} = O(1) \, $
as $ \,t\,\tend\,\infty $.
Now, given any sequence
$ \, t_{\mbox{}_{\scriptstyle n}} \,\tend\, \infty \, $
such that
$ \, {\xi}_{\mbox{}_{\mbox{}_{\scriptstyle n}}} \equiv
{\xi}_{\mbox{}_{\mbox{}_{\scriptstyle t_{\mbox{}_{n}} }}} \!$
converges,
say $ \, {\xi}_{\mbox{}_{\mbox{}_{\scriptstyle n}}} \,\tend\,
      {\xi}_{\mbox{}_{\mbox{}_{\scriptstyle \ast}}} \, $,
we \linebreak
\mbox{} \vspace{-0.4cm} \\
then have, from (45), (47), \\
\mbox{} \vspace{-0.40cm} \\
\mbox{} \hspace{0.30cm}
$ {\displaystyle
\sqrt{\,t_n\,} \; \|\: \omega(\cdot,t_{\mbox{}_{\scriptstyle n}}) \,\nLinfty
\,=\;
\factormm \:
\frac{1}{\, \mu \,-\, \sigma \, h \,
  \mbox{\sf erf}\,({\xi}_{\mbox{}_{\mbox{}_{\scriptstyle n}}}) \,}
\, \int_{0}^{1} \expxiyan \!\! \Hzeroy \; dy
} $, \\
\nl
\mbox{} \vspace{-0.30cm} \\
so that, letting $\, n \,\tend\,\infty $,
we obtain \\
\mbox{} \vspace{-0.40cm} \\
\mbox{} \hspace{0.20cm}
$ {\displaystyle
\lim_{\mbox{}_{\scriptstyle n\,\rightarrow\,\infty}}
\sqrt{\,t_n\,} \; \|\: \omega(\cdot,t_n) \,\nLinfty
\,=\;
\factormm \: \fbcm \:
\frac{ \mbox{\sf e}^{\mbox{}^{\scriptstyle
         -\, \mbox{\footnotesize $ {\xi}_{\mbox{}_{\scriptstyle \ast}}
                                        ^{\mbox{}^{\scriptstyle 2}} $} }} }
     {\, \mu \,-\, \sigma \, h \,
        \mbox{\sf erf}\,({\xi}_{\mbox{}_{\scriptstyle \ast}}) \,}
} $. \\
\nl
\mbox{} \vspace{-0.20cm} \\
This gives \\
\mbox{} \vspace{-0.40cm} \\
(48) \hspace{1.00cm}
$ {\displaystyle
\limitsup \;\; \thalf \: \|\: \omega(\cdot,t) \,\nLinfty
\,\leq\;
\factormm \: \fbcm \: \|\, {\cal F} \,\nLinfty
} $, \\
\nl
\mbox{} \vspace{-0.30cm} \\
which, together with (47) above, shows $(iii)$.
\hfill $\Box$ \\
\np
One consequence from Theorem 2 which is worth mentioning
it explicitly is the following one. \\
\nl
%
% ---------------------------------------------------------- Theorem 3
{\bf Theorem 3} \\
{\em
Let $\; a, \,b,\, c,\, \hat{a},\, \hat{b},\, \hat{c}\; $
be real constants, with $ \, c,\, \hat{c} > 0 $,
and let $\, u(x,t), \, \hat{u}(x,t) \,$ be the solutions
of equations $\:(17)$, $(18)$, respectively,
corresponding to initial states $ \, u_0, \, \hat{u}_0 \,$ in $\: \Lone \,$
with the same mass $\: m \neq 0$.
Then the following statements are equivalent to one another: \\
\mbox{} \vspace{-0.1cm} \\
$(i)$ \hspace{0.4cm}
    $ (\, a, \, b,\, c\,) \,=\, (\,\hat{a},\,\hat{b},\,\hat{c}\,) $, \\
\mbox{} \vspace{-0.05cm} \\
$(ii)$ \hspace{0.25cm}
$ {\displaystyle
\limitinf \;\; \tpp \: \|\, u(\cdot,t) \,-\, \hat{u}(\cdot,t) \,\nLp
\,=\; 0
} $ \hspace{0.30cm}
for some $ \, 1 \leq p \leq \infty $, \\
\mbox{} \vspace{-0.0cm} \\
$(iii)$ \hspace{0.20cm}
$ {\displaystyle
\limit \;\; \tpp \: \|\, u(\cdot,t) \,-\, \hat{u}(\cdot,t) \,\nLp
\,=\; 0
} $ \hspace{0.30cm}
for all \/ $ \, 1 \leq p \leq \infty $. \\
}
\nl
% -------------------------------------------------------------- Theorem 3
%
\mbox{} \vspace{-0.10cm} \\
$\Box$
Using Lemma 3, it is sufficient to examine the case
$\, u_0 = v_0 = m \, \mbox{\tiny $\mbox{}^{\mbox{\normalsize $\chi$}}$}
                     _{\mbox{}_{\scriptstyle [\,0,\,1\,]}} $,
$ m \neq 0 $.
If $ \, a \neq \hat{a} $, then,
from (6), (7), (8), there exist constants $\, K, \, k > 0 \,$
such that \\
\mbox{} \vspace{-0.1cm} \\
\mbox{} \hspace{3.25cm} $ {\displaystyle
|\; u(\, \xi\,\sqrt{\,t\,} \,+\, \hat{a}\,t\,,\, t\,) \;|
\;\leq\;
\frac{\,K\,}{\,\sqrt{\,t\,}\,} \:
\mbox{\sf e}^{\mbox{}^{\scriptstyle -\,
\frac{\mbox{$ \scriptstyle 1 $}}{\mbox{$ \scriptstyle 2 $}} \,
\mbox{\footnotesize $ (\, a \,-\, \hat{a} \,)^{\mbox{}^{\scriptstyle 2}}
                      \, t $} }}
                        } $ \\
\mbox{} \vspace{-0.2cm} \\
and \\
\mbox{} \vspace{-0.5cm} \\
\mbox{} \hspace{4.00cm} $ {\displaystyle
|\; \hat{u}(\, \xi\,\sqrt{\,t\,} \,+\, \hat{a}\,t\,,\, t\,) \;|
\;\geq\; \frac{\,k\,}{\,\sqrt{\,t\,}\;}
                        } $ \\
\nl
\mbox{} \vspace{-0.3cm} \\
for all $\, |\,\xi\,| \leq 1 \, $ and $ \, t \geq 1 $.
This clearly gives \\
\mbox{} \vspace{-0.1cm} \\
\mbox{} \hspace{4.0cm} $ {\displaystyle
\|\: u(\cdot,t) \,-\, \hat{u}(\cdot,t) \:\nLp
\;\geq\;
\kappa \, \tmp
                       } $ \\
\nl
for all $\,t\,$ large and $ \, 1 \leq p \leq \infty $,
for some constant $\, \kappa > 0 $. \\
\mbox{} \vspace{-0.25cm} \\
Assuming now that $\, a = \hat{a} $, suppose we have
$ \, (\, b,\,c \,) \neq (\,\hat{b},\,\hat{c}\,) \,$:
from (11), (16), we can find $ \, 1 < \pstar < \infty \,$
such that the limits (1) corresponding to $ \, u(x,t) \,$
and $ \,\hat{u}(x,t) $ \linebreak
are different, i.e.,
$ {\gamma}_{\mbox{}_{\mbox{}_{\scriptstyle p_{\mbox{}_{\mbox{}_{\!\ast}}} }}}
\!\neq\;
{\hat{\gamma}}_{\mbox{}_{\mbox{}_{\scriptstyle
  p_{\mbox{}_{\mbox{}_{\ast}}} }}} \!
$, where, for every $\,p$, \\
\mbox{} \vspace{-0.10cm} \\
\mbox{} \hspace{1.00cm}
$ {\displaystyle
\, {\gamma}_{\mbox{}_{\mbox{}_{\scriptstyle p }}} \,=\;
\limit \;\; \tpp \: \|\, u(\cdot,t) \,\nLp
, \hspace{0.75cm}
\, {\hat{\gamma}}_{\mbox{}_{\mbox{}_{\scriptstyle p }}} \,=\;
\limit \;\; \tpp \: \|\, \hat{u}(\cdot,t) \,\nLp
} $. \\
\nl
\mbox{} \vspace{-0.30cm} \\
In particular, we get \\
\mbox{} \vspace{-0.20cm} \\
(49a) \hspace{1.00cm}
$ {\displaystyle
\limitinf \;\; t^{\mbox{}^{\scriptstyle
                        \frac{\mbox{$ \scriptstyle 1 $}}
                             {\mbox{$ \scriptstyle 2 $}} \,
               \left(\, 1 \,-\, \frac{\mbox{$ \scriptstyle 1 $}}
                    {\mbox{$ \scriptstyle p_{\mbox{}_{\mbox{}_{\ast}}} $}}
               \right) }} \|\, u(\cdot,t) \,-\, \hat{u}(\cdot,t) \,
      \|_{\mbox{}_{\mbox{}_{\scriptstyle p_{\mbox{}_{\mbox{}_{\ast}}} }}}
\;\geq\;
|\, {\gamma}_{\mbox{}_{\mbox{}_{\scriptstyle p_{\mbox{}_{\mbox{}_{\ast}}} }}}
\,-\,
{\hat{\gamma}}_{\mbox{}_{\mbox{}_{\scriptstyle
  p_{\mbox{}_{\mbox{}_{\ast}}} }}} \,|
\,>\, 0
} $. \\
\nl
\mbox{} \vspace{-0.3cm} \\
Given $ \, p > \pstar $, we have, by interpolation of
$ \, \| \cdot \nLone $, $ \| \cdot \nLp $ at $ \,\pstar $, \\
\mbox{} \vspace{-0.20cm} \\
\mbox{} \hspace{0.40cm}
$ {\displaystyle
\|\, u(\cdot,t) \,-\, \hat{u}(\cdot,t)
\,\|_{\mbox{}_{\mbox{}_{\scriptstyle p_{\mbox{}_{\mbox{}_{\ast}}} }}}
\,\leq\;
\|\, u(\cdot,t) \,-\, \hat{u}(\cdot,t)
\,\|_{\mbox{}_{\mbox{}_{\scriptstyle 1}}}
    ^{\mbox{}^{\mbox{}^{\scriptstyle
      \frac{\,\mbox{$ \scriptstyle p \,-\, p_{\mbox{}_{\mbox{}_{\ast}}} $}\,}
           {\,\mbox{$ \scriptstyle p \,-\, 1 $}\,}
   \: \frac{\,\mbox{$ \scriptstyle 1 $}\,}
           {\,\mbox{$ \scriptstyle p_{\mbox{}_{\mbox{}_{\ast}}} $}} }}}
\:
\|\, u(\cdot,t) \,-\, \hat{u}(\cdot,t)
\,\|_{\mbox{}_{\mbox{}_{\scriptstyle p}}}
    ^{\mbox{}^{\mbox{}^{\scriptstyle
      \frac{\,\mbox{$ \scriptstyle p_{\mbox{}_{\mbox{}_{\ast}}} -\, 1 $}\,}
           {\,\mbox{$ \scriptstyle p \,-\, 1 $}\,}
   \: \frac{\,\mbox{$ \scriptstyle p $}\,}
           {\,\mbox{$ \scriptstyle p_{\mbox{}_{\mbox{}_{\ast}}} $}} }}}
} $ \\
\nl
\mbox{} \vspace{-0.40cm} \\
for every $ \, t > 0 $.
This gives, from (49a) above, \\
\np
\mbox{} \vspace{-0.70cm} \\
(49b) \hspace{0.75cm}
$ {\displaystyle
\limitinf \;\; \tpp \: \|\, u(\cdot,t) \,-\, \hat{u}(\cdot,t) \,\nLp
\;\geq\;
C \;
|\, {\gamma}_{\mbox{}_{\mbox{}_{\scriptstyle p_{\mbox{}_{\mbox{}_{\ast}}} }}}
\,-\,
{\hat{\gamma}}_{\mbox{}_{\mbox{}_{\scriptstyle
   p_{\mbox{}_{\mbox{}_{\ast}}} }}}
\, |^{\mbox{}^{\scriptstyle \mbox{\large $($}\, 1 \,-\,
     \frac{\,\mbox{$ \scriptstyle 1 $}\,}{\,\mbox{$ \scriptstyle p $}\,}
     \,\mbox{\large $)$} \:
     \frac{\,\mbox{$ \scriptstyle p_{\mbox{}_{\mbox{}_{\ast}}} $}\,}
          {\,\mbox{$ \scriptstyle p_{\mbox{}_{\mbox{}_{\ast}}} -\, 1 $}\,}
     }}
\!,
} $ \\
\nl
\mbox{} \vspace{-0.40cm} \\
where
$ {\displaystyle \; C \,=\,
     (\, {\gamma}_{\mbox{}_{\mbox{}_{\scriptstyle 1}}} +\,
         {\hat{\gamma}}_{\mbox{}_{\mbox{}_{\scriptstyle 1}}}
     \,)^{\mbox{}^{\scriptstyle -\,
     \mbox{\large $($} \, 1 \,-\,
         \frac{\,\mbox{$ \scriptstyle p_{\mbox{}_{\mbox{}_{\ast}}} $} }
              {\:\mbox{$ \scriptstyle p $} \,}
     \,\mbox{\large $)$} \:
         \frac{\,\mbox{$ \scriptstyle 1 $}\,}
         {\,\mbox{$ \scriptstyle p_{\mbox{}_{\mbox{}_{\ast}}} -\, 1 $}\,}
     }} \!
} $.
Similarly, for $ \, p < \pstar $, we get \\
\mbox{} \vspace{-0.10cm} \\
\mbox{} \hspace{1.25cm}
$ {\displaystyle
\|\, u(\cdot,t) \,-\, \hat{u}(\cdot,t)
\,\|_{\mbox{}_{\mbox{}_{\scriptstyle p_{\mbox{}_{\mbox{}_{\ast}}} }}}
\,\leq\;
\|\, u(\cdot,t) \,-\, \hat{u}(\cdot,t)
\,\|_{\mbox{}_{\mbox{}_{\scriptstyle p}}}
    ^{\mbox{}^{\mbox{}^{\scriptstyle
      \frac{\,\mbox{$ \scriptstyle p $}\,}
           {\,\mbox{$ \scriptstyle p_{\mbox{}_{\mbox{}_{\ast}}} $}} }}}
\:
\|\, u(\cdot,t) \,-\, \hat{u}(\cdot,t)
\,\|_{\mbox{}_{\mbox{}_{\scriptstyle \infty}}}
    ^{\mbox{}^{\mbox{}^{\scriptstyle  1 \,-\,
      \frac{\,\mbox{$ \scriptstyle p $}\,}
           {\,\mbox{$ \scriptstyle p_{\mbox{}_{\mbox{}_{\ast}}} $}} }}}
\!,
} $ \\
\nl
\mbox{} \vspace{-0.40cm} \\
which gives, using (49a), \\
\mbox{} \vspace{-0.0cm} \\
(49c) \hspace{0.35cm}
$ {\displaystyle
\limitinf \;\; \tpp \: \|\, u(\cdot,t) \,-\, \hat{u}(\cdot,t) \,\nLp
\;\geq\;
|\, {\gamma}_{\mbox{}_{\mbox{}_{\scriptstyle p_{\mbox{}_{\mbox{}_{\ast}}} }}}
\,-\,
{\hat{\gamma}}_{\mbox{}_{\mbox{}_{\scriptstyle
   p_{\mbox{}_{\mbox{}_{\ast}}} }}}
\, |^{\mbox{}^{\scriptstyle
     \frac{\,\mbox{$ \scriptstyle p_{\mbox{}_{\mbox{}_{\ast}}} $}}
          {\:\mbox{$ \scriptstyle p $}\,} }}
\:
(\, {\gamma}_{\mbox{}_{\mbox{}_{\scriptstyle \infty}}} +\,
    {\hat{\gamma}}_{\mbox{}_{\mbox{}_{\scriptstyle \infty}}}
\,)^{\mbox{}^{\scriptstyle -\,
         \frac{\,\mbox{$ \scriptstyle p_{\mbox{}_{\mbox{}_{\ast}}} -\, p $}}
              {\:\mbox{$ \scriptstyle p $} \,} }}
\!.
} $ \\
\nl
\mbox{} \vspace{-0.30cm} \\
Hence, in all cases above,
$ \, (\, a,\, b,\, c \,) \neq (\,\hat{a},\, \hat{b},\, \hat{c} \,) \,$
gives, for every $ \, 1 \leq p \leq \infty $, \\
\mbox{} \vspace{-0.10cm} \\
\mbox{} \hspace{3.50cm}

$ {\displaystyle
\limitinf \;\; \tpp \: \|\, u(\cdot,t) \,-\, \hat{u}(\cdot,t) \,\nLp
\; > \; 0
} $, \\
\nl
\mbox{} \vspace{-0.30cm} \\
which, together with Theorem 2, finishes the argument.
\hfill $\Box$ \\
%
% -------------------------------------------- END of proof for Theorem 3
\nl
In a similar way, we can show that, given initial states
$ \, u_0 , \, \tilde{u}_0 \in \Lone \, $
with different masses,
the corresponding solutions $ \, u(\cdot,t), \, \tilde{u}(\cdot,t) \,$
of equation (2) satisfy, for every $ \, 1 \leq p \leq \infty $, \\
\mbox{} \vspace{-0.20cm} \\
(50) \hspace{3.50cm}
$ {\displaystyle
\|\, u(\cdot,t) \,-\, \tilde{u}(\cdot,t) \,\nLp
\;\geq\; c_p \, \tmp
} $ \\
\nl
for all $ \, t \, $ large, where $ \, c_p $ is some positive
constant. \\
\nl
% -----------------------------------------------------
\nl
%
% ----------------------------------------------------------------------
%
%                            REFERENCES
%
% ----------------------------------------------------------------------
%
\nl
{\large {\bf References}}  \\
\nl
\mbox{} \vspace{-0.3cm} \\
\mbox{[\hspace{1mm}1\hspace{1mm}]} \hspace{1mm}
            \begin{minipage}[t]{14.0cm}
              E. R. Benton and G. W. Platzman,
              {\em A table of solutions of the one-dimensional
              Burgers equation},
              Quart. Appl. Math. {\bf 30} (1972), 195 -- 212.
            \end{minipage} \\
\nl
\mbox{} \vspace{-0.3cm} \\
\mbox{[\hspace{1mm}2\hspace{1mm}]} \hspace{1mm}
            \begin{minipage}[t]{14.0cm}
              Z. Brze\'zniak and B. Szafirski,
              {\em Asymptotic behaviour of $ L^{\mbox{}^{\scriptstyle 1}}\!$
              norm of solutions to parabolic equations},
              Bull. Polish Acad. Sci. Math. {\bf 39} (1991), 1 -- 10.
            \end{minipage} \\
\nl
\mbox{} \vspace{-0.3cm} \\
\mbox{[\hspace{1mm}3\hspace{1mm}]} \hspace{1mm}
            \begin{minipage}[t]{14.0cm}
              J. M. Burgers, {\em Application of a model system
              to illustrate some points of the statistical
              theory of free turbulence},
              Nederl. Akad. Wefensh. Proc. {\bf 43} (1940), 2 -- 12.
%             Proc. Acad. Sci., Amsterdam, {\bf 43} (1940), 2 -- 12.
            \end{minipage} \\
\nl
\mbox{} \vspace{-0.3cm} \\
\mbox{[\hspace{1mm}4\hspace{1mm}]} \hspace{1mm}
            \begin{minipage}[t]{14.0cm}
              J. D. Cole, {\em On a quasilinear parabolic equation
              occurring in aerodynamics},
              Quart. Appl. Math. {\bf 9} (1951), 225 -- 236.
            \end{minipage} \\
\nl
\mbox{} \vspace{-0.3cm} \\
\mbox{[\hspace{1mm}5\hspace{1mm}]} \hspace{1mm}
            \begin{minipage}[t]{14.0cm}
              E. Hopf, {\em The partial differential equation
              $ \,u_{t} + u\,u_{x} = \mu\,u_{xx} $},
              Comm. Pure Appl. Math. {\bf 3} (1950), 201 -- 230.
            \end{minipage} \\
\nl
\mbox{} \vspace{-0.3cm} \\
\mbox{[\hspace{1mm}6\hspace{1mm}]} \hspace{1mm}
            \begin{minipage}[t]{14.0cm}
              R. Rudnicki, {\em Asymptotic stability in
              $L^{\mbox{}^{\scriptstyle 1}}\!$ of $\,$parabolic
              equations},
              J. Diff. Equations, {\bf 102} (1993), 391 -- 401.
            \end{minipage} \\
\nl
\mbox{} \vspace{-0.3cm} \\

\end{document}